\documentclass[11pt]{article}

\usepackage[margin=2.5cm]{geometry}
\usepackage{amsmath, amsthm, amsfonts, amssymb}
\usepackage{graphicx, color, xcolor}
\usepackage{hyperref}
\usepackage[utf8]{inputenc}
\usepackage{cleveref}
\usepackage[sort&compress,numbers,square]{natbib}
\usepackage{algorithm} 
\usepackage{algpseudocode} 
\usepackage{mathrsfs} 
\usepackage{orcidlink} 
\usepackage{caption} 
\usepackage{subcaption} 
\usepackage{url} 
\usepackage{booktabs} 

\title{\sc{Adaptive Sensor Placement Inspired by Bee Foraging: Towards Efficient Environment Monitoring} \\[1ex] \large}

\author{\Large Sai Krishna Reddy Sathi\orcidlink{0009-0005-4136-3005}\thanks{Indian Institute of Technology Madras, Chennai - 600036, India. \texttt{me21b181@smail.iitm.ac.in}}}

\date{}

\begin{document}
\maketitle

\hrule 
\begin{abstract}
This paper aims to make a mark in the future of sustainable robotics, where efficient algorithms are required to carry out tasks like environmental monitoring and precision agriculture efficiently. We proposed a hybrid algorithm that combines Artificial Bee Colony (ABC) with Levy flight to optimize adaptive sensor placement alongside an important notion of hotspots from domain knowledge experts. By enhancing exploration and exploitation, our approach significantly improves the identification of critical hotspots. This versatile algorithm also holds promise for broader search and rescue operations applications, demonstrating its potential in optimization problems across various domains.

\vspace{5pt}
		
\end{abstract}

\vspace{10pt}
\hrule
	
\section{Introduction}
\label{sec:intro}
The challenges of climate change and environmental degradation call for innovative solutions that can support sustainable practices across various fields, including robotics. Sustainable robotics aims to create robotic systems that can operate in harmony with nature, focusing on designs that reduce environmental impact and support ecosystem health. The role of sustainable robotics in the different aspects of acheiving sustainability was researched extensively by Bugmann et al. \cite{bugmann2011role}. Environmental monitoring is one of the crucial domains among the many areas in sustainable robotics. Forest canopies, in particular, are critical for biodiversity, water cycles, and carbon storage. Monitoring these areas helps scientists understand ecosystem dynamics, assess forest health, and identify threats, yet accessing them is difficult without causing disruption. This is where sustainable, biomimetic robotic systems can make a difference, providing efficient, minimally invasive solutions for complex environmental monitoring tasks.

Biomimetic design—drawing inspiration from nature—plays a significant role in sustainable robotics. Researchers like Prof. Mirko Kovac and Prof. Stefano Mintchev have advanced the use of nature-inspired models in this field. For instance, Prof. Kovac’s work at Imperial College focuses on creating adaptive robotic systems that mimic animal behaviors, allowing robots to operate effectively in natural settings. Prof. Mintchev similarly investigates bio-inspired strategies that enable robots to adapt to their surroundings. These principles help create efficient and environmentally friendly robots, making them suitable for long-term monitoring of delicate ecosystems like rainforests. Kirchgeorg et al. developed robotic systems like AVOCADO \cite{kirchgeorg2023design} and probes to collect temperature, humidity data, bioacoustic signatures, eDNA, etc, for forest canopy monitoring in the environmental robotics laboratory and similarly in the SMA-based dart-like sensor attachment mechanisms are developed in Prof. Kovac's lab for environmental monitoring \cite{hamaza2020sensor}. These systems are of primary relevance to this research and are assumed to be individual robotic units in this swarm for covering the forest.

\begin{figure}[H]
    \centering
    \begin{subfigure}[b]{0.45\textwidth}
        \centering
        \includegraphics[width=\textwidth]{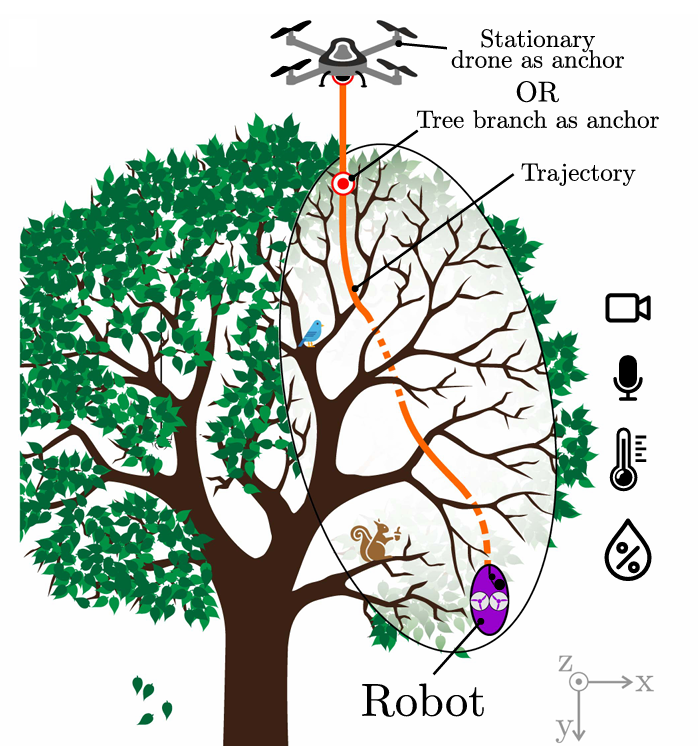} 
        \caption{Approach of the AVOCADO robot with a drone as an anchor, picture is taken from \cite{kirchgeorg2023design}}
        \label{fig:figure1}
    \end{subfigure}
    \hfill
    \begin{subfigure}[b]{0.45\textwidth}
        \centering
        \includegraphics[width=\textwidth]{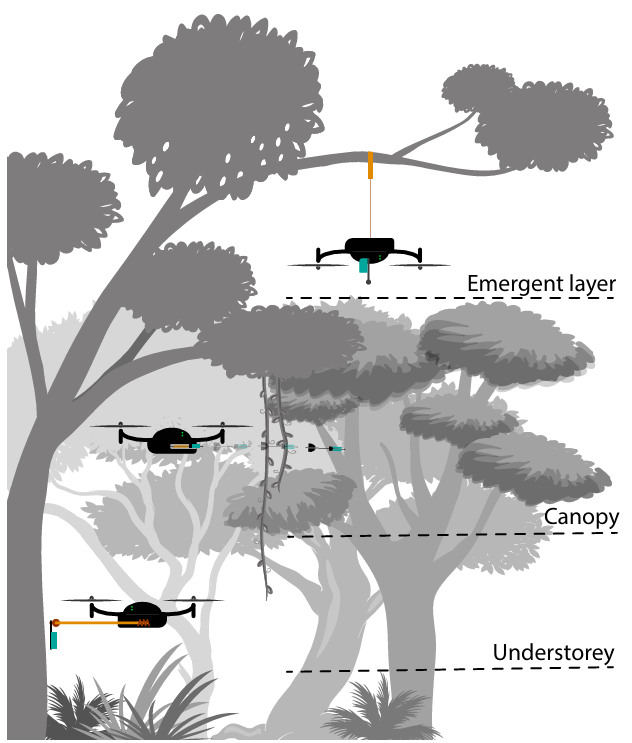} 
        \caption{Approach for sensor mounting on the tree from Prof. Kovac's lab, picture is taken from \cite{hamaza2020sensor}}
        \label{fig:figure2}
    \end{subfigure}
    \label{fig:two_figures}
\end{figure}

A key part of this biomimetic approach involves adopting strategies from nature for navigating and covering complex areas. Bee-foraging behaviors, for example, have inspired several algorithms for optimizing search and resource-gathering tasks. In nature, bees efficiently explore large areas and communicate routes to each other to balance exploration of new areas with thorough searches of known ones. The Artificial Bee Colony (ABC) algorithm, which models this behavior, divides tasks between scout bees (exploring) and forager bees (collecting information). This balance of exploration and focused search makes ABC-based methods particularly useful for applications like sensor placement, where coverage is essential, but energy must be conserved.

\section{Problem Statement}
This paper presents a new application of bee-foraging algorithms for adaptive sensor placement in forest canopy exploration. We use a hybrid Adaptive Bee Colony-Levy algorithm that combines ABC’s search dynamics with the Levy flight method, which introduces random but efficient search patterns. Our goal is to optimize sensor distribution to detect as many areas of interest as possible while avoiding redundancy by collecting data around hotspot locations in the forest whose locations are roughly estimated by the domain experts (environmental researchers in this case). By applying these principles, this work contributes to sustainable robotics and lays the groundwork for conservation, resource management, and climate monitoring applications.

For a standard reference, we chose the XPrize Rainforest competition. In the final, the teams are supposed to survey 100 hectares of rainforest in 24 hrs. Active research has been going on building robotic systems for sensor deployment and collecting data in these dense canopies \cite{kirchgeorg2023design} \cite{hamaza2020sensor}. In this work, we consider a swarm of autonomous versions of these robots and figure out the optimal approach for efficient data collection using a hotspot approach. A hotspot, as used in this paper's context, is an area within the forest canopy that shows high biodiversity or unusual ecological activity. By employing the hotspot approach, we collect high-quality data where we can extract crucial insights from the forest. Since the data regarding the time taken for this robotic unit to deploy a sensor module in one location and other parameters regarding the UAV, in this paper, the absolute values of the time taken for covering are not relevant. The aim is to cover the forest grid, prioritizing the area around hotspots, which is effectively an optimization problem. In the simulation, the forest grid is considered to be a 100x100 cell grid with all the UAVs initialized from the center of the bottom edge of the forest grid, i.e., at (50,0).  

\begin{figure}[H]
    \centering
    \includegraphics[width=\textwidth]{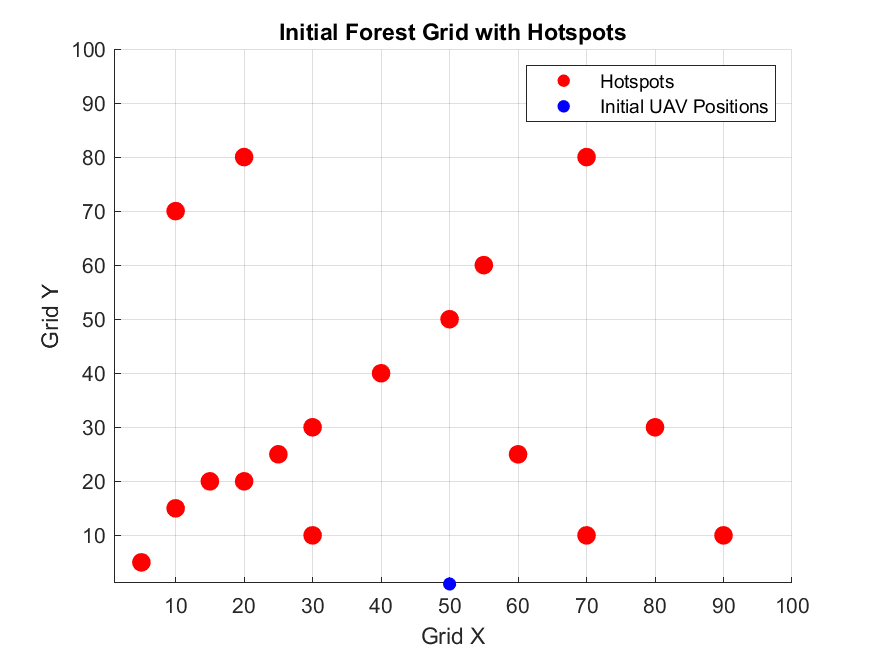}  
    \caption{Forest-Grid initialization with 20 hotspots with all the UAVs at (50,0)}
    \label{fig:unet_architecture}
\end{figure}

\section{Literature Review}

Several algorithms were researched for solving the MOOP (Multi-Objective Optimization Problem) of Multi-UAV path planning. Spanning several fields, evolutionary algorithms like GA (Genetic Algorithm), PCO (Particle Swarm Optimization), and ACO (Ant-Colony Algorithm) have been widely used. Bee-foraging algorithms have proven effective in a variety of environmental exploration applications due to their inherent balance of exploration and exploitation. The Artificial Bee Colony (ABC) algorithm, initially proposed by Karaboga (2007), simulates the foraging behavior of bees, using a combination of employed, onlooker, and scout bees to locate and refine search areas \cite{karaboga2007powerful}. Given a population of $N$ bees, the search space is defined as $X = \{ x_1, x_2, \ldots, x_N \}$, with each bee exploring a region based on its nectar probability $P_i$ defined as:
\begin{equation}
    P_i = \frac{fit_i}{\sum_{j=1}^{N} fit_j},
\end{equation}
where $fit_i$ represents the fitness of solution $i$. This probabilistic model enhances the ability of robots to monitor environmental features effectively.

Many studies focus only on improving local search or strengthening its global search capabilities in a singular manner, without considering a balanced approach to both \cite{wiesemuller2021zero}. Chiang et al. enhanced the ABC algorithm by combining it with discretized honey source optimization and support vector machine parameters to improve classification accuracy and convergence speed. However, this approach reduces the local search ability \cite{chiang2020novel}. Strategies like those proposed by Alatas, which incorporate neighborhood search and Chaos theory, aim to strengthen local search but do not address the global search imbalance \cite{alatas2010chaotic}. Deniz et al. combined ABC with Differential Evolution to reduce the risk of local optima during early searches \cite{ustun2017study}, while Li et al. used a roulette-based strategy to enhance the algorithm’s ability to escape local minima \cite{li2020orthodontic}. 

As expected, our initial attempts at solving the problem using standard ABC or PSO methods failed by getting stuck in a local optimum in certain configurations of the hotspots. Hence, we proceeded to research the implementation of a hybrid approach with Levy flights.
Several studies have successfully tackled the challenge of standard Artificial Bee Colony (ABC) algorithms getting trapped in local optima by integrating Levy flights, which introduce a stochastic, long-step search component to improve global exploration. Yang and Deb \cite{yang2009cuckoo} demonstrated that Levy flights, with their unique long-tailed distributions, enhance global search effectiveness by preventing premature convergence, especially in complex landscapes. Liu et al. \cite{liu2018artificial} proposed a dynamic penalty-based ABC algorithm with Levy flights, showing that the hybrid approach successfully navigates constrained optimization problems by balancing exploration and exploitation. Based on these successful applications, leveraging and tuning a hybrid ABC-Levy approach for our problem is well-justified, as it enhances search dynamics and addresses local optima challenges effectively.

Through hybridizing ABC with Levy flight models. Levy flights, characterized by a probability distribution with heavy tails, facilitate long-distance exploration. The hybrid ABC-Levy model enables a bee to execute random walks defined by:
\begin{equation}
    x_{t+1} = x_t + \alpha L(s, \lambda),
\end{equation}
where $L(s, \lambda) \sim |s|^{-(1+\lambda)}$ for step size $s$ and scale parameter $\alpha$, achieving a blend of local and global search capability 

Due to standard ABC and PSO algorithm limitations, our hybrid model was designed to improve adaptability and reliability in dense forests. By introducing a weighted balance factor $\omega$ between exploration and exploitation, our approach builds on these foundational works to achieve optimal environmental monitoring. The Levy-enhanced ABC offers a promising solution for hotspot identification, both in forest canopies and broader environmental monitoring applications.



\subsection{Methodology}

In designing our approach, we adopted a methodology that balances localized intensive monitoring with broader, exploratory scanning. We define key terms here that will serve as the foundation for the framework.

\begin{itemize}
    \item \textbf{Hotspot:} A specific region identified within the forest canopy that exhibits heightened ecological significance, such as high biodiversity or specific environmental features like rare plant clusters or animal nesting zones. This parameter is central to prioritizing areas for in-depth study. Further research has to be done towards building an error-free hotspot identification system from a practical implementation perspective. One plausible approach can be having a camera mounted on the UAV, pointing downwards (towards the canopy), and developing deep learning-based models. For instance, we could be looking for canopy-based estimations similar to the work done by Prof. Kovoc but from the top view of canopies.  They could also be any other ecological markers from visual feed that are of importance, which environmental researchers can identify.
    
    \item \textbf{Scouting:} In the context of our framework, scouting refers to the preliminary search by robotic agents for potential hotspots. Inspired by bee-foraging patterns, scouting is conducted by robotic "scouts" that identify regions with high probability of ecological interest based on sensor readings, often through random or probabilistic paths. 
    
\end{itemize}

Our methodology employs a hybrid approach of bee-inspired adaptive sensing, blending random exploration with hotspot-focused monitoring. Our robotic agents execute a two-phase operation: \textit{scouting} phase, where initial data is gathered over a large area, and \textit{focused monitoring} phase, where hotspots are revisited to gather detailed measurements.

Let $H = \{h_1, h_2, \ldots, h_n\}$ represent identified hotspots in the canopy, and let $s_i$ denote scouting paths. Each robot is assigned a function $\phi(h_i, s_i)$ that evaluates the richness of biodiversity indicators in a region. The overall goal is to maximize the cumulative biodiversity metric $B$, calculated as:
\begin{equation}
    B = \sum_{i=1}^{n} \phi(h_i, s_i),
\end{equation}
where $\phi$ represents the diversity or density of detected species.

Our adaptive methodology utilizes a bee-foraging inspired algorithm, which allows each robot to dynamically adjust its search path based on real-time data inputs. Such an approach provides flexibility to cover large areas while still dedicating resources to promising regions. 


\section{Bee-Foraging Algorithms}

\subsection{Artificial Bee Colony (ABC) Algorithm}
The Artificial Bee Colony (ABC) algorithm divides bees into three categories: employed bees, onlooker bees, and scout bees. Let \(\vec{x}_{i} \in \mathbb{R}^{n}\) denote the position of the \(i\)-th bee in an \(n\)-dimensional search space, where the objective function \(f(\vec{x})\) represents the fitness of each solution.

\begin{equation}
    f(\vec{x}_{i}) = \text{objective function value of position } \vec{x}_{i}
\end{equation}

Each bee generates a new candidate position \(\vec{v}_{i}\) around its current position \(\vec{x}_{i}\):
\begin{equation}
    \vec{v}_{i} = \vec{x}_{i} + \phi_{i,j} (\vec{x}_{i} - \vec{x}_{k})
\end{equation}
where \(\phi_{i,j} \sim \mathcal{U}(-1, 1)\) is a uniform random number, and \(k \neq i\) is a randomly selected solution index.

After evaluating \(f(\vec{v}_{i})\), the new position is accepted if it yields a better fitness:
\begin{equation}
    \text{if } f(\vec{v}_{i}) < f(\vec{x}_{i}), \text{ then } \vec{x}_{i} = \vec{v}_{i}.
\end{equation}

\subsection{Particle Swarm Optimization (PSO)}
The Particle Swarm Optimization (PSO) algorithm updates each particle’s position \(\vec{x}_{i}\) and velocity \(\vec{v}_{i}\) based on its personal best position \(\vec{p}_{i}\) and the global best position \(\vec{g}\). The velocity update equation is given by:
\begin{equation}
    \vec{v}_{i}(t+1) = \omega \vec{v}_{i}(t) + c_{1} r_{1} (\vec{p}_{i} - \vec{x}_{i}(t)) + c_{2} r_{2} (\vec{g} - \vec{x}_{i}(t)),
\end{equation}
where \(\omega\) is the inertia weight, \(c_{1}\) and \(c_{2}\) are acceleration coefficients, and \(r_{1}, r_{2} \sim \mathcal{U}(0,1)\) are random factors.

The position update is:
\begin{equation}
    \vec{x}_{i}(t+1) = \vec{x}_{i}(t) + \vec{v}_{i}(t+1).
\end{equation}

The above algorithms were tested in the forest grid map with randomly initialized hotspots. It was observed that they were stuck at local minima points in several cases, leading to inefficient search. For example, consider 20 hotspots and 10 of them are closer to the side of the grid the UAVs start from, and the other 10 are almost as a cluster on the other side of the forest grid, then it was seen that our swarm of UAVs almost just get stuck exploring the closer 10 hotspots and regions around them. Hence, it was decided to add an exploration component to the exploitation component. In other words, it is asked to take a random walk associated with a weight factor with the traditional approach. The hybrid of Levy Flight with PSO has shown poor performance despite tuning several hyperparameters. Therefore, we proceeded further with the Hybrid ABC-Levy flight algorithm.

\subsection{Hybrid ABC with Levy Flight}
In the hybrid ABC approach, we incorporate a Levy flight mechanism to improve exploration. A Levy flight is defined by a step length \(L\) drawn from a Levy distribution:
\begin{equation}
    L \sim \text{Levy}(\alpha, \beta) \text{ where } 0 < \alpha \leq 2 \text{ and } \beta = 1.
\end{equation}

The step for each bee is updated as:
\begin{equation}
    \vec{x}_{i} = \vec{x}_{i} + L \cdot (\vec{x}_{i} - \vec{x}_{k}),
\end{equation}
where \(L\) ensures a heavy-tailed distribution allowing for both local and global search.

The probability of selecting Levy flight in hybrid ABC is given by a control parameter \(\lambda\), which balances exploration and exploitation:
\begin{equation}
    \lambda = \frac{1}{1 + e^{-\sigma (f(\vec{x}_{i}) - f(\vec{g}))}}
\end{equation}
where \(\sigma\) adjusts sensitivity to differences between current and global best fitness values.

This approach improves the convergence rate by balancing exploration (Levy flight) with exploitation (ABC updates), reducing the chance of local minima trapping.

\begin{algorithm}[H]
\caption{Hybrid ABC-Levy Flight Algorithm for Adaptive Sensor Placement}
\label{alg:abc-levy}
\begin{algorithmic}[1]
\Require Environment grid $\mathcal{G}$ with hotspot set $\mathcal{H}$, exploration and exploitation parameter $\alpha$, exploitation parameter $\beta$, safe-zone radius $r$, Levy weight $\lambda$

\State \textbf{Initialize} UAVs at initial positions $\{x_i\}_{i=1}^n$ within $\mathcal{G}$, set parameters $\alpha$, $\beta$, $r$, $\lambda$

\While{termination criteria not met (e.g., all hotspots covered)}
    \State \textbf{Employed Bee Phase: Foraging and Exploration}
    \State 1. \textbf{Levy Flight Movement:} 
    \State \hspace{1em} Each UAV \( x_i \) moves based on the Levy distribution:
    \vspace{-2pt}
    \[
    L_i = \lambda \cdot u \cdot |v|^{-1/\beta}, \quad u, v \sim N(0, 1)
    \]
    \vspace{-1pt}
    where $u, v$ are samples from a Gaussian distribution, generating a random movement vector with long-tail distribution, promoting exploration over long distances.

    \State 2. \textbf{Position Update via Levy Step:}
    \State \hspace{1em} Update each UAV’s position:
   \[
    x_i \leftarrow x_i + L_i
    \]
    \State \hspace{1em} Apply boundary conditions to constrain $x_i$ within $\mathcal{G}$:
    \[
    x_i \leftarrow \max(\min(x_i, \text{max boundary}), \text{min boundary})
    \]
    
    \State 3. \textbf{Exploration-Exploitation Balancing:}
    \State \hspace{1em} For high-fitness UAVs, increase local search focus:
    \vspace{-2pt}
    \[
    x_i \leftarrow x_i + \beta \cdot (x_i - x_j), \quad x_j \text{ is nearby high-fitness UAV}
    \]
    \hspace{2.5em} For lower-fitness UAVs, encourage exploration:
    \vspace{-2pt}
    \[
    x_i \leftarrow x_i + \alpha \cdot (x_{\text{best}} - x_i), \quad x_{\text{best}} \text{ is the best-known position}
    \]

    \State 4. \textbf{Safe-Zone Constraint Enforcement:}
    \State \hspace{1em} Ensure UAVs avoid clustering within $r$ distance:
    \[
    x_i \leftarrow x_i + r \cdot \frac{x_i - x_j}{\|x_i - x_j\|}, \quad \|x_i - x_j\| < r
    \]

    \State \textbf{Onlooker Bee Phase: Selective Intensification}
    \State 5. \textbf{Reinforced Levy Update for Select Onlookers:}
    \State \hspace{1em} Onlooker UAVs are selected with probability proportional to their coverage performance. 
    \hspace{5em} They execute Levy-based exploration:
    \vspace{-2pt}
    \[
    x_i^{\text{new}} = x_i + \lambda \cdot u \cdot |v|^{-1/\beta}
    \]

    \State 6. \textbf{Evaluate Fitness of New Position:}
    \State \hspace{1em} Fitness is calculated based on proximity to unvisited hotspots:
    \[
    f(x_i^{\text{new}}) = \sum_{k \in \mathcal{H}} w_k \cdot \mathbb{I}(\|x_i^{\text{new}} - h_k\| \leq r), \quad w_k \text{ is the importance of hotspot } h_k
    \]
    Update $x_i \leftarrow x_i^{\text{new}}$ if fitness improves.

    \State \textbf{Scout Phase:} If coverage stagnates, reset $x_i$ to a random position to avoid local maxima.
\EndWhile

\State \Return Optimized sensor positions $\{x_i\}_{i=1}^n$
\end{algorithmic}
\end{algorithm}

\section{Results and Discussion}
The initial simulation results from the standard PSO and ABC very quickly proved to be ineffective upon changing the hotspot configurations. Hence, further discussion will be on the results from the Hybrid ABC-Levy algorithm. The additional characteristics that were imparted in the code are as follows:

\begin{itemize}

     \item \textbf{Collision Avoidance:} During this process, a collision of any two or more UAVs would be totally undesirable. In our decentralized swarm, the concept of an artificial potential field is implemented to ensure no collisions despite the commands given by the ABC-Levy algorithm.
     
    \item \textbf{Avoiding no hotspot zones:} Since the locations of hotspots are assumed to be already known, if the drone follows our hybrid ABC-Levy algorithm and enters a zone where there are no hotspots within the threshold radius (manually decided depending on the hotspot distribution), it immediately tries to get out of the circle, moving towards the nearest hotspot.

     \item \textbf{Adjusting the simulation to match realistic drone dynamics:} This is effectively adjusting the distance that can be traversed by a node (UAV) for each step update. Manually, it was observed that a MaxStepSize of 5 units seemed reasonable. More accurate results can be obtained, given the dynamics model of each UAV. 
\end{itemize}

In the demonstration of the results below, 20 hotspots were considered with 5 UAVs. One of the main parameters to be tuned is the LevyWeight. The behavior of the swarm upon changing LevyWeight has been carefully observed and manually tuned, alongside the time taken for the coverage metric. The best results were observed for LevyWeight around 3 for different hotspot configuration types. A dynamic LevyWeight value based on the hotspot configuration is left to future works.

\begin{figure}[H]
    \centering
    \resizebox{1\textwidth}{!}{ 
        \begin{subfigure}[b]{0.45\textwidth}
            \centering
            \includegraphics[width=\textwidth]{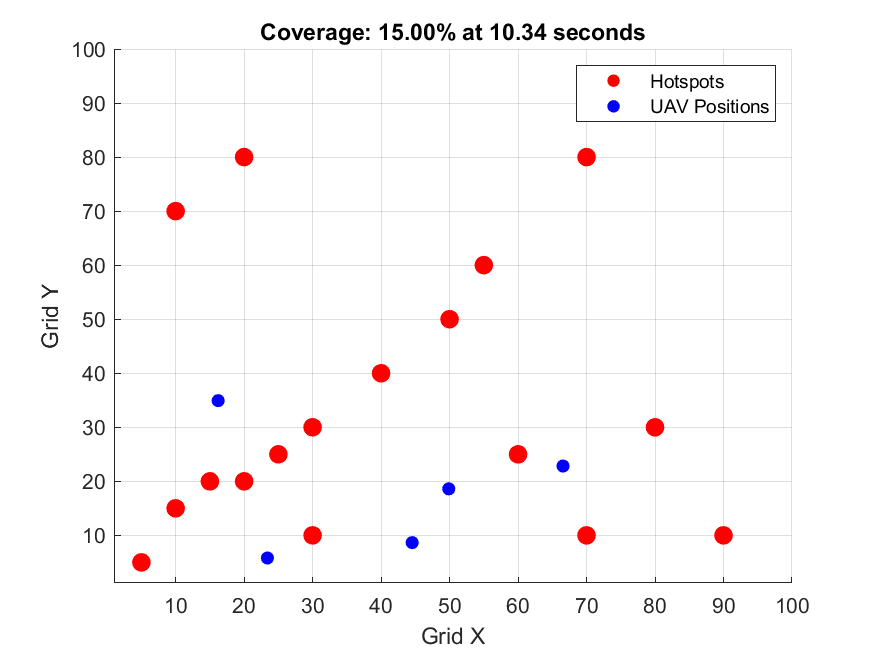}
            \label{fig:subfig1}
        \end{subfigure}
        \hfill
        \begin{subfigure}[b]{0.45\textwidth}
            \centering
            \includegraphics[width=\textwidth]{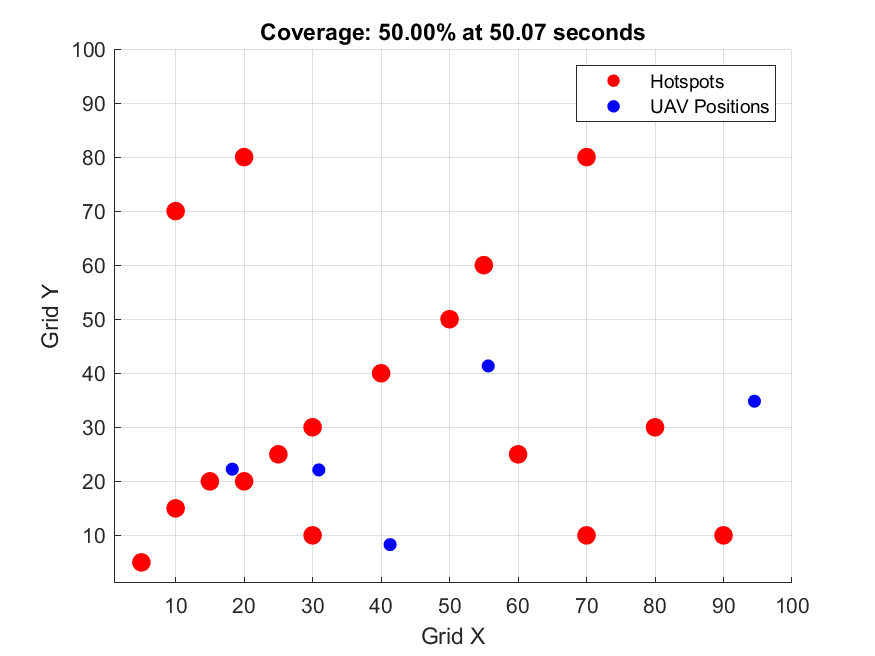}
            \label{fig:subfig2}
        \end{subfigure}

        \vspace{1em} 

        \begin{subfigure}[b]{0.45\textwidth}
            \centering
            \includegraphics[width=\textwidth]{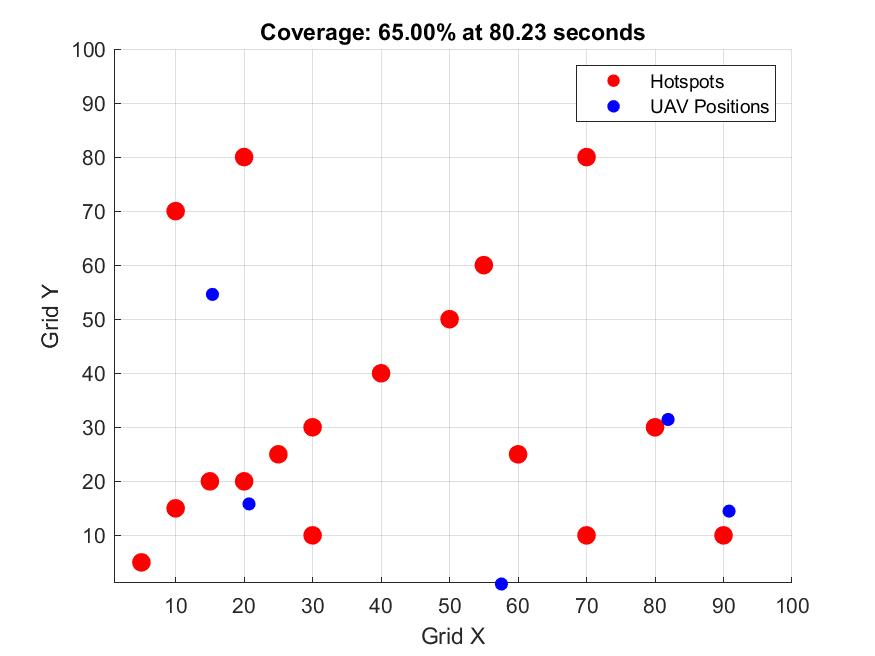}
            \label{fig:subfig3}
        \end{subfigure}
        \hfill
        \begin{subfigure}[b]{0.45\textwidth}
            \centering
            \includegraphics[width=\textwidth]{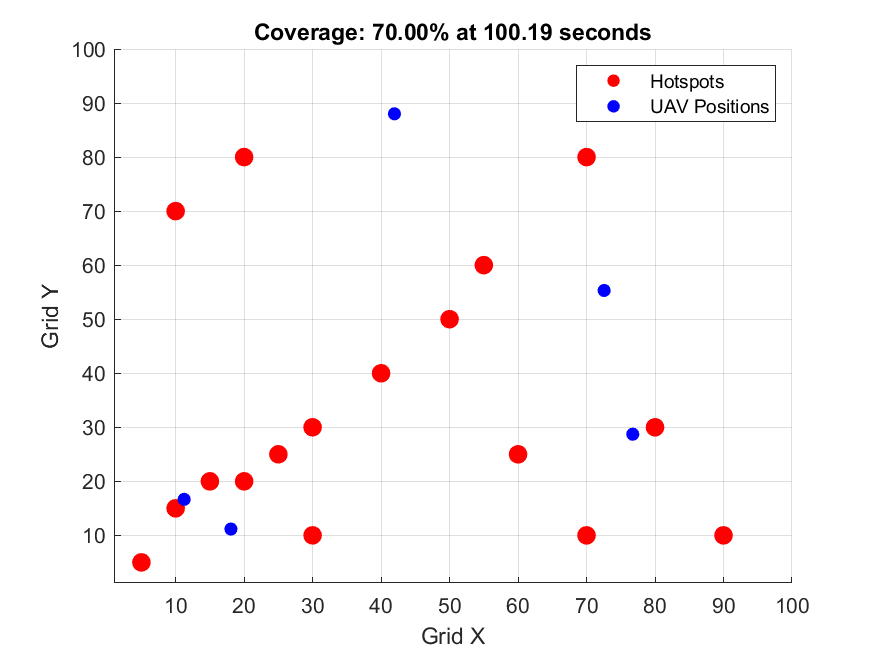}
            \label{fig:subfig4}
        \end{subfigure}
    }
    
    \caption{Frames of the simulation video with LevyWeight = 3; \textbf{Time taken to cover all hotspots: 227.5006 seconds}}
    \label{fig:main_figure}
\end{figure}

\begin{figure}[H]
    \centering
    \resizebox{1\textwidth}{!}{ 
        \begin{subfigure}[b]{0.45\textwidth}
            \centering
            \includegraphics[width=\textwidth]{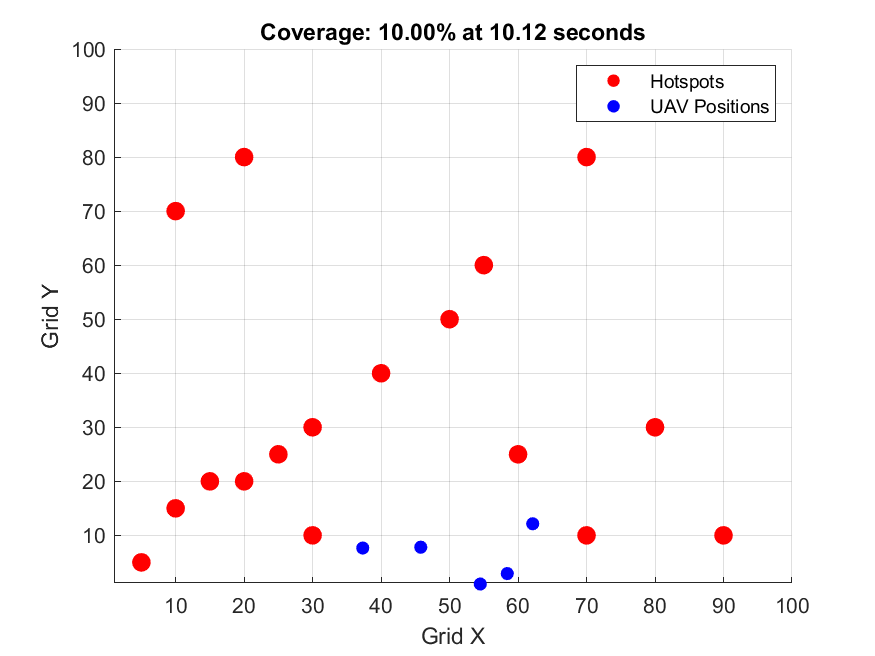}
            \label{fig:subfig1}
        \end{subfigure}
        \hfill
        \begin{subfigure}[b]{0.45\textwidth}
            \centering
            \includegraphics[width=\textwidth]{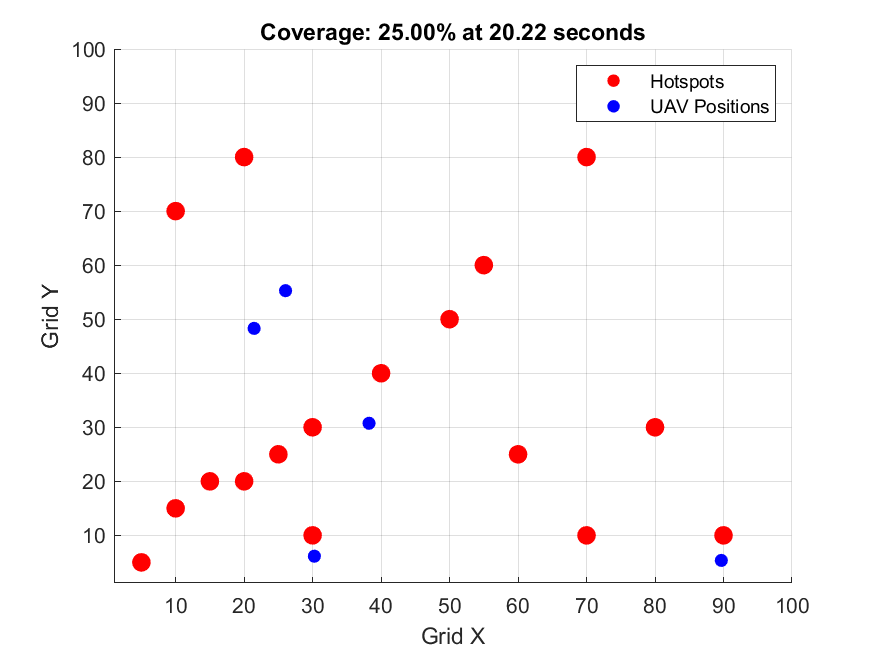}
            \label{fig:subfig2}
        \end{subfigure}

        \vspace{1em} 

        \begin{subfigure}[b]{0.45\textwidth}
            \centering
            \includegraphics[width=\textwidth]{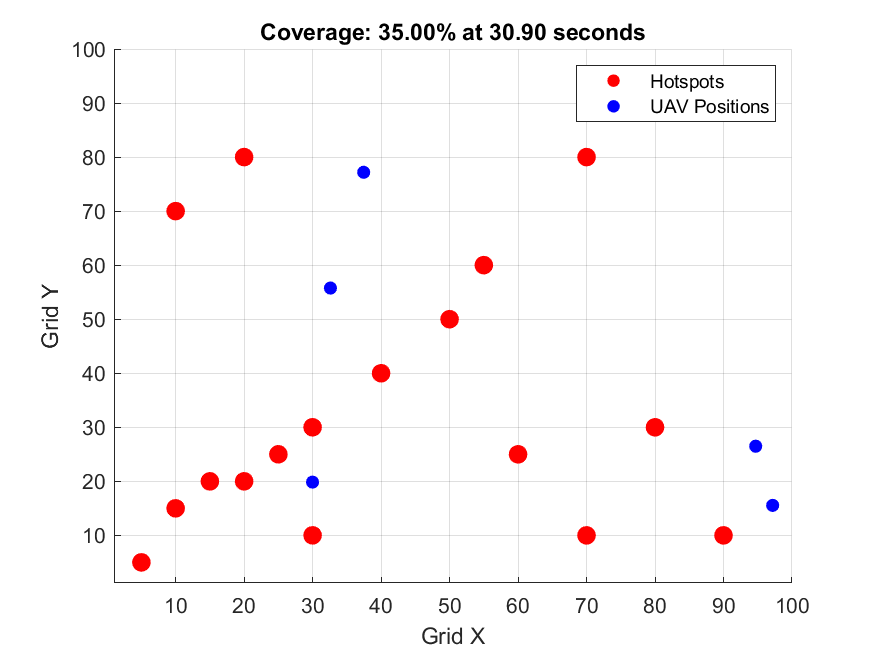}
            \label{fig:subfig3}
        \end{subfigure}
        \hfill
        \begin{subfigure}[b]{0.45\textwidth}
            \centering
            \includegraphics[width=\textwidth]{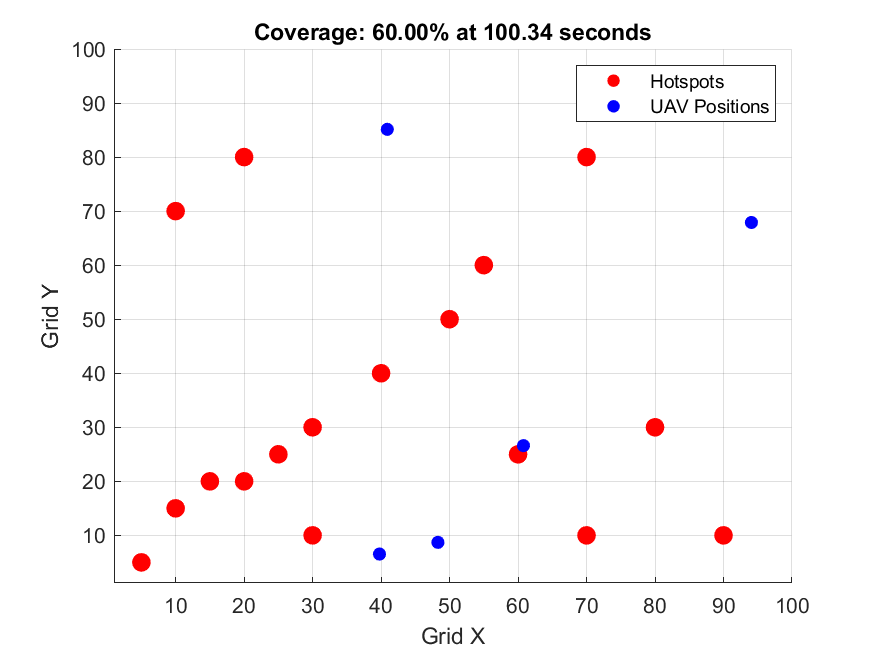}
            \label{fig:subfig4}
        \end{subfigure}
    }
    
    \caption{Frames of the simulation video with LevyWeight = 5; \textbf{Time taken to cover all hotspots: 350.9228 seconds}}
    \label{fig:main_figure}
\end{figure}

A heatmap is also plotted for different LevyWeight values, it effectively conveys the grid cell coverage density across the forest grid. The case following the hotspot spatial distribution pattern would be the best fit. The heatmap metric gives the same optimal value of 3, complimenting the time taken metric shown previously. Heatmaps were analyzed for higher weight values as well, which tend to be ineffective because of excessive random walks overshadowing the ABC algorithm.

\begin{figure}[H]
    \centering
    \includegraphics[width=\textwidth]{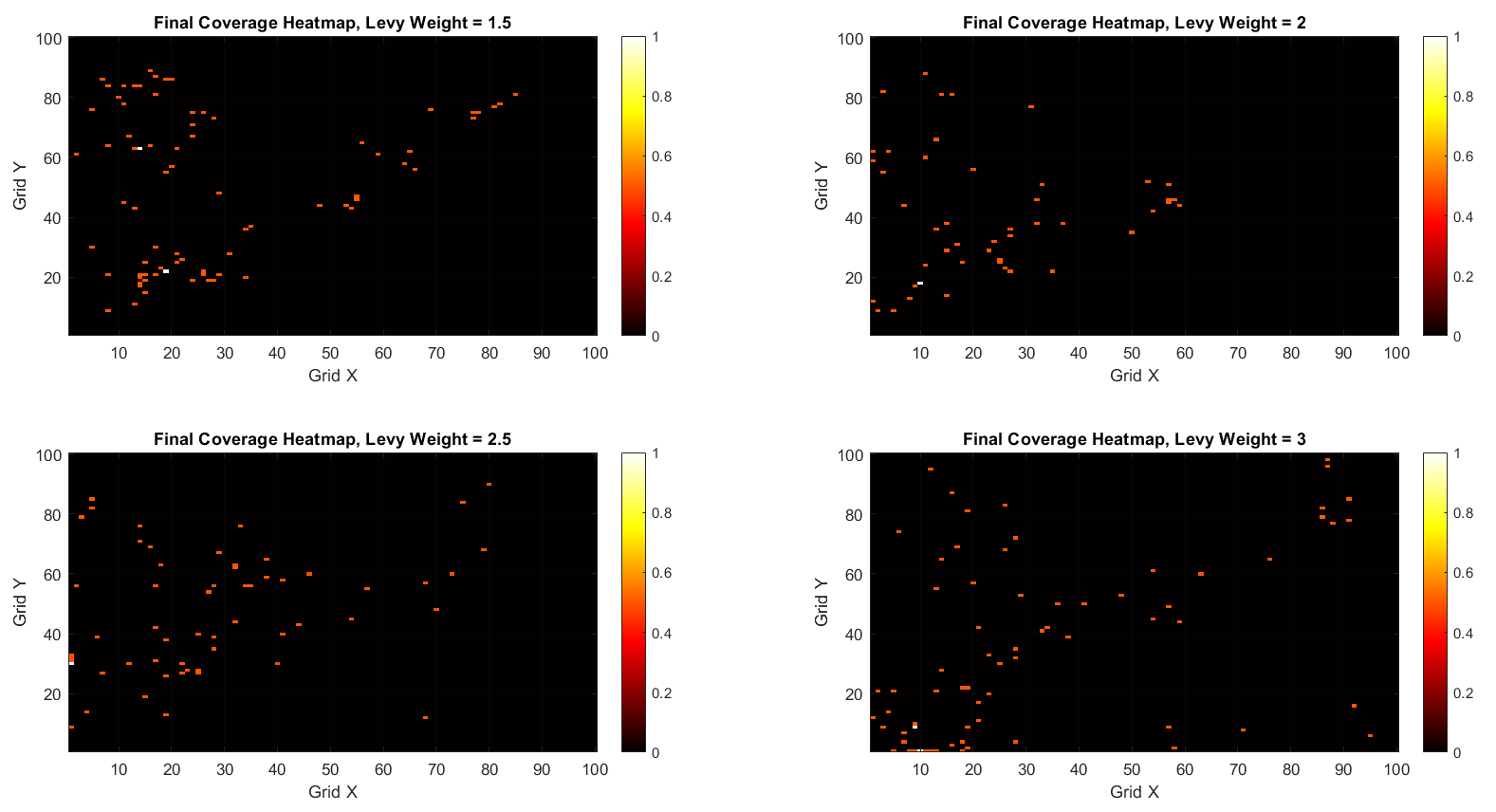}  
    \caption{Coverage Heatmap for LevyWeights 1.5, 2, 2.5 and 3 respectively}
    \label{fig:unet_architecture}
\end{figure}


\section{Conclusion and Future Works}
This paper presented a hybrid Artificial Bee Colony (ABC) and Levy flight algorithm to address the challenge of optimizing sensor placement within a simulated environment with hotspots. The results demonstrate promising capabilities in achieving efficient sensor coverage by balancing exploration and exploitation. This adaptive algorithm showcases considerable robustness in simulated conditions; however, deploying it in real-world environments presents further challenges, especially in GPS-denied, dense forest settings where localization and mapping remain difficult.

Localization and mapping are critical for fully autonomous robotic systems navigating these environments. Since GPS signals can be unreliable in forests, integrating algorithms like this with SLAM (Simultaneous Localization and Mapping) technologies may improve feasibility. SLAM enables robots to create maps in real-time while tracking their position within it, a vital feature for our approach in unknown terrains. Prior works highlight how advancements in visual SLAM and LiDAR-based mapping could assist our system in dense forests by providing more precise position estimates despite the absence of GPS.

Detecting hotspots in dense, GPS-denied environments requires innovative strategies. One effective approach is multi-sensor fusion, integrating data from thermal, infrared, and audio sensors to enhance environmental perception. Additionally, techniques like visual odometry and LiDAR mapping can facilitate accurate localization and mapping.

Swarm intelligence algorithms offer another promising method, enabling multiple robots to collaborate in systematically exploring an area and sharing hotspot data. Furthermore, employing machine learning can allow the system to adapt its strategies based on previously encountered hotspots, improving detection efficiency over time.

An adaptive control system could optimize sensor placements based on environmental feedback, such as temperature or sound variations, to enhance real-time hotspot detection. While these strategies show potential, their practical effectiveness will require thorough field testing and iterative adjustments to ensure reliability in real-world applications

This algorithm also holds potential beyond environmental monitoring. Applications such as search and rescue, where targeted detection and monitoring of specific areas are essential, could benefit from this approach. Our optimization framework could adapt dynamically in these high-stakes scenarios, ensuring critical areas are rapidly identified and communicated to rescue teams. Further, the algorithm's adaptability and hotspot-detection capabilities align well with disaster response and hazard mapping.

In future works, we aim to enhance the algorithm's real-time adaptability by exploring deeper integration with dynamic SLAM methods and by evaluating its performance on autonomous aerial platforms equipped with sensors that can operate under dense canopy cover. Another direction involves developing an error-correction layer for instances where autonomous platforms may lose track of mapped hotspots. By refining these features, our hybrid algorithm could advance toward fully autonomous field deployments, enabling real-world applications across diverse and challenging environments.

\newpage
\bibliography{mybib}

\end{document}